\documentclass[11pt]{elsarticle}

\usepackage{enumerate}
\usepackage{graphicx}
\usepackage{float}
\usepackage{hyperref}
\usepackage[margin=1in]{geometry}
\usepackage{comment}
\usepackage{color}

\usepackage{amssymb, color}
\usepackage{amsmath}
 \usepackage{pifont}
 \usepackage[utf8]{inputenc}

 \usepackage[utf8]{inputenc}
\usepackage[T1]{fontenc}
\usepackage{enumitem}

\makeatletter
\def\LaTeX{\leavevmode L\raise.42ex
    \hbox{\kern-.3em\size{\sf@size}{0pt}\selectfont A}\kern-.15em\TeX}
\makeatother

\newcommand{\BibTeX}{{\rm B\kern-.05em{\sc
          i\kern-.025emb}\kern-.08em\TeX}}

\makeatletter
\def\@currentlabel{2.1}\label{e:dispaa}
\def\@currentlabel{2.21}\label{e:dispau}
\def\@currentlabel{2.22}\label{e:dispav}
\def\@currentlabel{2.23}\label{e:dispaw}
\def\@currentlabel{2.24}\label{e:dispax}
\def\theequation{\thesection.\@arabic\c@equation}
\makeatother

\renewcommand{\theequation}{\arabic{section}.\arabic{equation}}
\def \D{\Delta}

\newcommand{\N}{\mathbb N}
\newcommand{\e}{\epsilon}

\def \D{\Delta}

\newtheorem{thm}{Theorem} [section]
\newtheorem{lem}{Lemma} [section]
\newtheorem{prop}{Proposition} [section]

\newtheorem{rem}{Remark}[section]

\renewcommand{\theequation}{\thesection.\arabic{equation}}
\renewcommand{\thesection}{\arabic{section}}
\renewcommand{\theequation}{\thesection.\arabic{equation}}
\let\ssection=\section\renewcommand{\section}{\setcounter{equation}{0}\ssection}

\begin{document}

\begin{frontmatter}

\title{On the regularity and partial regularity of extremal solutions of a Lane-Emden system}

\author[]{Hatem Hajlaoui}
\ead{hajlouihatem@gmail.com}
\address{Institut Pr\'{e}paratoire aux \'Etudes d'Ing\'{e}nieurs. Universit\'e de Kairouan, Tunisie.}
\begin{abstract} In this paper,  we consider the system $-\Delta u =\lambda (v+1)^p,\;\;-\Delta v = \gamma (u+1)^\theta$ on a smooth
bounded domain $\Omega$ in $\mathbb{R}^N$ with the Dirichlet boundary condition $u=v=0$ on $\partial \Omega.$ Here $ \lambda,\gamma$ are positive parameters.
Let $x_0$ be the largest root of the polynomial
\begin{equation*}
H(x) = x^4 -
\frac{16p\theta(p+1)(\theta+1)}{(p\theta-1)^2}x^2 +
\frac{16p\theta(p+1)(\theta+1)(p+\theta+2)}{(p\theta-1)^3}x
-\frac{16p\theta(p+1)^2(\theta+1)^2}{(p\theta-1)^4}.
\end{equation*} We show that the extremal solutions associated to the above system are bounded provided
  $N<2+2x_0.$ This improves the previous work in \cite{co1}. We also prove that, if $N\geq 2+2x_0,$ then the singular set of any extremal
  solution has Hausdorff dimension less or equal to $N-(2+2x_0).$
\end{abstract}
\begin{keyword}
Extremal solution\sep stable solution\sep regularity and partial regularity.
\end{keyword}
\end{frontmatter}
\section{Introduction}
\setcounter{equation}{0}

In this paper, we are interested in the regularity and partial regularity of extremal solutions to the following system:
 \begin{eqnarray}\label{1.1}
  \left\{ \begin{array}{rll}
   -\Delta u  &= \lambda (v+1)^p & \hbox{in } \Omega  \\
 -\Delta v &= \gamma (u+1)^\theta & \hbox{in }\Omega  \\
u &= v =0 & \hbox{on } \partial\Omega,
\end{array}\right.
  \end{eqnarray}
  where $\Omega$ is a smooth bounded domain in $\mathbb{R}^N$, $\lambda, \gamma$ are positive parameters and $p,\theta > 1$. In particular, we examine when the extremal solutions  of \eqref{1.1} are smooth. Applying standard elliptic theory, it is sufficient to show that the extremal solutions are bounded.
 The nonlinearities we examine naturally fit into the following general assumptions:
$$\mbox{$f$ is smooth, positive, increasing, convex in $[0, \infty)$, and }\lim_{u \rightarrow \infty} \frac{f(u)}{u} =\infty.
\leqno{({\rm R})}$$

Recalling that the scalar analog of the system \eqref{1.1} is given by
$$
-\Delta u = \lambda f(u) \;\; \hbox{in }\Omega, \quad
u = 0 \;\; \hbox{on } \partial\Omega. \leqno{(Q)_\lambda}
$$
When the nonlinearity $ f$ satisfies (R), it
is well known that there exists a finite positive critical parameter $ \lambda^*$ such that
for all $ 0< \lambda < \lambda^*$ there exists a smooth minimal solution $ u_\lambda$ of $ (Q)_\lambda$. Here the minimal solution means in the pointwise sense. In addition, the minimal solution $ u_\lambda$ is semi-stable in the sense that
\[ \int_\Omega \lambda f'(u_\lambda) \psi^2 dx \le \int_\Omega | \nabla \psi|^2 dx, \quad \forall\; \psi \in H_0^1(\Omega).\]
Moreover, the map $ \lambda \mapsto u_\lambda(x)$ is increasing on $[0,\lambda^*).$ This allows one to define $ u^*(x):= \lim_{\lambda \nearrow \lambda^*}
u_\lambda(x)$, the so-called extremal solution, which is shown to be the unique weak solution of $(Q)_{\lambda^*}$, and there is no weak solution of
$(Q)_\lambda$ for $ \lambda > \lambda^*$. The regularity and properties of extremal solution to $(Q)_\lambda$ have attracted a lot of attention. It is known that it depends on the nonlinearity $f$,
the dimension $N$ and the geometry of the domain $\Omega$. See for instance \cite{bcmr,BV,Cabre,CC,Martel,MP,Nedev, v, yz}.

\medskip
The situation is much less understood for the corresponding elliptic system. Consider the generalization of
\eqref{1.1} as follows:
$$
\left\{
\begin{array}{rll}
-\Delta u &= \lambda f(v)&  \hbox{in }\Omega  \\
-\Delta v &= \gamma g(u)  & \hbox{in }\Omega  \\
u &= v =0 & \hbox{on } \partial\Omega.
\end{array}\right.
\leqno{(P)_{\lambda, \gamma}}
$$
 Define  $ \mathcal{Q} :=\{ (\lambda,\gamma): \lambda, \gamma >0 \}$ and
\[ \mathcal{U}:= \left\{ (\lambda,\gamma) \in \mathcal{Q}: \mbox{ there exists a smooth solution $(u,v)$ of $(P)_{\lambda,\gamma}$} \right\}.\]
Set $\Upsilon:= \partial \mathcal{U} \cap \mathcal{Q}$, which plays the role of the extremal parameter $ \lambda^*$. As shown by Montenegro  \cite{Mont},
if $f$ and $g$ satisfy (R), then
\begin{enumerate}
\item $\mathcal{U}$ is nonempty. For all $ (\lambda,\gamma) \in \mathcal{U}$, there is a minimal solution of $(P)_{\lambda,\gamma}$.
\item For each $ 0 < \sigma < \infty$ there is some $ 0 < \lambda^*_\sigma < \infty$ such that   $ \mathcal{U} \cap \{ (\lambda,\sigma \lambda): 0 < \lambda \}$  is given by  $ \{ (\lambda, \sigma \lambda): 0 < \lambda < \lambda_\sigma^* \} \cup \mathcal{H}$
 where $\mathcal{H}$ is either the empty set or $ \{ (\lambda_\sigma^*, \sigma \lambda_\sigma^*) \}$. The map $\sigma \mapsto \lambda_\sigma^*$ is bounded
 on compact subsets of $(0,\infty)$.
 \item Fix $ 0 < \sigma < \infty$ and let $(u_\lambda,v_\lambda)$ denote the minimal solution of $(P)_{\lambda,\sigma \lambda}$ for
 $0 <\lambda < \lambda_\sigma^*$.  Then $u_\lambda,v_\lambda$ are increasing in $ \lambda$ and
\begin{align}
\label{new0}
u^*(x):= \lim_{\lambda \nearrow \lambda_\sigma^*} u_\lambda(x), \quad v^*(x):= \lim_{\lambda \nearrow \lambda_\sigma^*} v_\lambda(x)
\end{align}
 is always a weak solution to $(P)_{\lambda_\sigma^*, \sigma \lambda_\sigma^*}$.
\end{enumerate}
 In addition, let $ (\lambda,\gamma) \in \mathcal{U}, $ the minimal solution $(u,v)$ of $(P)_{\lambda,\gamma}$ is semi-stable in the sense that there are
 $0 < \zeta,\chi \in H_0^1(\Omega)$ and $\eta \ge 0$  such that
\begin{equation} \label{1.2}
 -\Delta \zeta = \lambda  f'(v) \chi + \eta \zeta, \;\; -\Delta \chi = \gamma g'(u) \zeta+ \eta \chi \quad  \mbox{ in } \Omega.
 \end{equation}
See \cite{Mont} and also \cite{co1} for an alternative proof of \eqref{1.2}. Moreover, we have the following useful inequality, see Lemma 1 in \cite{co1}
and Lemma 3 in \cite{dfs}.
\begin{lem} \label{stabb}
 Let $(u,v)$ denote a semi-stable solution of $(P)_{\lambda,\gamma}$ in the sense of (\ref{1.2}). Then
\begin{equation} \label{1.3}
 \sqrt{\lambda \gamma} \int_\Omega \sqrt{f'(v) g'(u)} \phi^2 \le \int_\Omega | \nabla \phi|^2, \quad \forall\; \phi \in H_0^1(\Omega).
\end{equation}
\end{lem}

For example, when $f(t) = g(t) = e^t$, it was shown in \cite{dfs} that for $1\leq N\leq 9$, the extremal solution $(u^*, v^*)$ is smooth, see also \cite{craig0, DG}. Furthermore,
if $N \geq 10$, D\'avila and Goubet showed that the Hausdorff dimension of the singular set of any extremal solution is less or equal to $N - 10$.
For the polynomial system \eqref{1.1}, Cowan proved in \cite{co1}:

\begingroup \renewcommand\thethm{\Alph{thm}}
\begin{thm}
\label{Cow}
Suppose that $ 1 < p \le \theta$, $(\lambda^*,\gamma^*) \in \Upsilon$. Then, the extremal solution $(u^*,v^*)$ of \eqref{1.1} is bounded provided
$N<2+\frac{4(\theta+1)t_0}{p\theta-1},$ where
    \begin{align}\label{1.4}
t_{0} = \sqrt{\frac{p\theta(p+1)}{\theta+1}}+\sqrt{\frac{p\theta(p+1)}{\theta+1}-\sqrt{\frac{p\theta(p+1)}{\theta+1}}}.
\end{align}
Consequently, the extremal solutions are smooth for any $1<p\leq \theta$ provided $N\leq 10.$
\end{thm} \endgroup

A main idea in \cite{co1} is to use the stability inequality \eqref{1.3}. This technique was used to consider various Liouville theorem and regularity of extremal
solutions for elliptic systems and biharmonic equations, see for example \cite{CEG, CG, co, hhm, hhy, dfs, craig0, craig2}.

\medskip
Our main concern here is to improve Cowan's result.
\setcounter{thm}{0}
\begin{thm} \label{main1} Let $(\lambda^*,\gamma^*) \in \Upsilon$ and $(u^*,v^*)$ denote the associated extremal solution of \eqref{1.1}. Suppose that
$N<2+2x_0,$ where $x_0$ be the largest root of the polynomial $H(x) = $
\begin{align}\label{1.5}
x^4 - \frac{16p\theta(p+1)(\theta+1)}{(p\theta-1)^2}x^2 +
\frac{16p\theta(p+1)(\theta+1)(p+\theta+2)}{(p\theta-1)^3}x
-\frac{16p\theta(p+1)^2(\theta+1)^2}{(p\theta-1)^4}.
\end{align}  Then $ u^*,v^*$ are bounded. In particular, the extremal solutions are smooth provided $N\leq 10.$
\end{thm}
Using Remark \ref{r.2.1} below, we see that $2t_0\frac{\theta+1}{p\theta-1}\leq x_0$ for any $1<p \leq \theta$, with equality if and only if $p=\theta$, where
  $t_0$ is given by \eqref{1.4}, so our result improves Theorem \ref{Cow}.

\medskip
To prove Theorem \ref{main1}, we will use the following Souplet type pointwise estimate between $u$ and $v$, solution of \eqref{1.1}. See Lemma 2 in \cite{co1}.
\begin{lem} \label{pointwise} Let $ (u,v)$ denote a smooth solution of \eqref{1.1} and suppose  that $ \theta \ge p >1$. Let
\[\alpha:= \max \left\{ 0, \left(  \frac{\gamma (p+1)}{\lambda (\theta+1) } \right)^\frac{1}{p+1}-1 \right\} .\] Then
\begin{equation} \label{point}
(v+1 +\alpha)^{p+1} \ge \frac{\gamma(p+1)}{\lambda(\theta +1)} (u+1)^{ \theta +1} \quad \mbox{ in $ \Omega$}.
\end{equation}
\end{lem}
Obviously, as $v > 0$ and $\alpha \geq 0$, we have
\begin{align}
\label{new1}
(v+1)^{p+1} \geq \left(\frac{v+1 +\alpha}{\alpha+1}\right)^{p+1} \geq \frac{\gamma(p+1)}{\lambda(\theta +1)(\alpha+1)^{p+1}} (u+1)^{ \theta +1} \quad \mbox{ in $ \Omega$}.
\end{align}

In the spirit of \cite{DG}, we are also interested in the partial regularity for extremal solutions. Let $(u^*, v^*)$ be an extremal solution of \eqref{1.1}, a
point $x \in \Omega$ is said regular if there exists a neighborhood of $x$ on which $u^{\star}$ and $v^{\star}$ are bounded; Otherwise $x$ is said singular. Denote by $\mathcal{S}$ the set of singular points of $(u^*,v^*).$ By definition, the regular set $\Omega\backslash{\mathcal S}$ is open and by elliptic regularity, $u^{\star}, \,v^{\star}$ are smooth in $\Omega\backslash{\mathcal S}$.
\begin{thm} \label{main2} Assume that $N\geq 2+2x_0,$ where $x_0$ is that in Theorem \ref{main1}. Let $(u^*,v^*)$ denote an extremal solution of \eqref{1.1}, i.e.~with $(\lambda^*, \gamma^*) \in \Upsilon$, then the Hausdorff dimension of its singular set $\mathcal{S}$ is less or
equal to $N-(2+2x_0).$
\end{thm}

\begin{rem} If $p=1$ or $\theta = 1$ and $p\theta > 1$, following the proofs of Theorems \ref{main1}-\ref{main2}, we can show that the results remain true. In other words, Theorems \ref{main1}-\ref{main2} hold true for $p, \theta \geq 1$ verifying $p\theta > 1$.
\end{rem}

This paper is organized as follows. We prove Theorem \ref{main1} in Section 2. The Section 3 is devoted to the proof of Theorem \ref{main2}.

\section{Proof of Theorem \ref{main1}}
First we remark that the polynomial $H$ is completely symmetric in $p$ and $\theta$. Hence we assume from now on $1 < p \leq \theta$, without loss of generality.

\medskip
The following Lemma plays an important role in dealing with Theorem \ref{main1}, where we use some ideas from \cite{hhm, hhy, co}. Let $(\lambda^*,\gamma^*)
\in \Upsilon$ and $ \sigma:= \frac{\gamma^*}{\lambda^*}$. Define $\Gamma_\sigma:=\{ (\lambda,\sigma \lambda): \frac{\lambda^*}{2} < \lambda < \lambda^* \}$ and denote $(u^*,v^*)$ the extremal solution associated to $(\lambda^*,\gamma^*)$ for \eqref{1.1} defined by \eqref{new0}, i.e.~$(u^*,v^*)$ is the pointwise limit of the minimal solutions along the ray $\Gamma_\sigma$ as $\lambda \nearrow \lambda^*$.
\begin{lem} \label{f}
 Let $ (\lambda^*,\gamma^*) \in \Upsilon$, $\sigma =\frac{\gamma^*}{\lambda^*}$, let $ (u,v)$ denote the minimal solution of \eqref{1.1} for $(\lambda, \gamma)\in \Gamma_\sigma$. Define
\begin{align}\label{2.1}
 L(s):=s^4-\frac{16p\theta(p+1)}{\theta+1}s^2+\frac{16p\theta(p+1)(p+\theta+2)}{(\theta+1)^2}s-\frac{16p\theta(p+1)^2}{(\theta+1)^2}.
\end{align}
Then for any $s > p+1$ verifying
$L(s) < 0$, there exists $C_s <\infty$ such that for any $(\lambda, \gamma)\in \Gamma_\sigma$, there holds
\begin{equation} \label{2.2}
\int_\Omega (u+1)^\frac{\theta-1}{2}(v+1)^\frac{p+2s-1}{2} + \int_\Omega (u+1)^{ \theta +  \frac{ (\theta+1) (s-1)}{p+1}} \le {C}_s.
\end{equation}
 \end{lem}

 \noindent
    \textbf{Proof.} We handle only the first integral in \eqref{2.2}, since the second estimate is an immediate consequence of the first one thanks to
    \eqref{new1}. Inserting $\phi:= (u+1)^{\frac{q+1}{2}}-1$ with $q > 0$ into (\ref{1.3}), we obtain
\begin{align}
\label{new3}
\sqrt{\lambda \gamma p \theta} \int_\Omega (v+1)^\frac{p-1}{2} (u+1)^\frac{\theta-1}{2} \left[(u+1)^\frac{q+1}{2}-1\right]^2 \le \frac{(q+1)^2}{4} \int_\Omega (u+1)^{q-1} | \nabla u|^2.
\end{align} On the other hand, multiplying the first equation of \eqref{1.1} by $(u+1)^{q}-1$, we get
\begin{align}
\label{new4} \frac{(q+1)^2}{4}\int_\Omega (u+1)^{q-1} | \nabla u|^2 = \frac{(q+1)^2\lambda}{4q}\int_\Omega (v+1)^p\Big[(u+1)^{q}-1\Big].\end{align}
Combining \eqref{new3}, \eqref{new4} and dropping some positive terms, there holds
\begin{align} \label{2.4}
\sqrt{ \lambda \gamma}a_1J_1 \leq \lambda \int_\Omega (v+1)^p (u+1)^{q} + 2\sqrt{\lambda \gamma}a_1 I_1,
\end{align}
 where
 $$a_1=\frac{4q\sqrt{p\theta}}{(q+1)^2}, \quad J_1 :=  \int_\Omega (v+1)^\frac{p-1}{2} (u+1)^\frac{\theta + 2q +1}{2},
 \quad I_1:= \int_\Omega (v+1)^\frac{p-1}{2} (u+1)^\frac{\theta+q}{2}.$$
 Similarly, using $\phi:= (v+1)^{\frac{r+1}{2}}-1$ in (\ref{1.3}) with $r > 0$, we obtain
\begin{align} \label{2.5}
\sqrt{ \lambda \gamma}a_2J_2 \le \gamma \int_\Omega (u+1)^\theta (v+1)^{r} + 2\sqrt{ \lambda \gamma}a_2 I_2
\end{align}
where
$$a_2=\frac{4r\sqrt{p\theta}}{(r+1)^2}, \quad J_2 := \int_\Omega (v+1)^\frac{p+2r+1}{2} (u+1)^\frac{\theta-1}{2}, \quad I_2:=\int_\Omega (v+1)^\frac{p+r}{2} (u+1)^\frac{\theta-1}{2}.$$
Fix now
\begin{equation}\label{2.6}
 q=\frac{(\theta+1)r}{p+1}+\frac{\theta-p}{p+1}, \quad \mbox{ or equivalently }
q+1=\frac{(\theta+1)(r+1)}{p+1}.
\end{equation}
Let $r> p$ and so $q>\theta,$ we claim that for any $\e > 0$, there exists $C_\e > 0$ independent of $(\lambda, \gamma) \in \Gamma_\sigma$ such that
\begin{eqnarray} \label{2.7}
I_1 \le \e J_1 + \e\int_\Omega (v+1)^p (u+1)^{q} + C_\e, \quad I_2 \le \e J_2 + \e \int_\Omega (u+1)^\theta (v+1)^{r} + C_\e.
\end{eqnarray}
Indeed, using successively Young's inequality for $(u+1)^\frac{\theta+q}{2}$ and $(v+1)^\frac{p-1}{2}$, we get
\begin{align*}
 I_1 \leq \e J_1 + C_\e \int_\Omega (v+1)^\frac{p-1}{2} & \leq \e J_1 + \e \int_\Omega (v+1)^p + C_\e\\
 & \le \e J_1 + \e\int_\Omega (v+1)^p (u+1)^{q} + C_\e.
\end{align*}
The estimate for $I_2$ is similar, so we omit it. Inserting \eqref{2.7} into \eqref{2.4} and \eqref{2.5} respectively, we get (for $\e < 1/2$)
\begin{equation*}
J_1 \leq \frac{1}{A_1}\int_\Omega (v+1)^p (u+1)^{q} + C_\e, \quad J_2 \leq \frac{1}{A_2}\int_\Omega (u+1)^\theta (v+1)^{r} + C_\e,
\end{equation*}
with
\[ A_1 = \frac{a_1(1 - 2\e)}{\sqrt{\frac{\lambda}{\gamma}} + 2a_1\e},\quad
 A_2 = \frac{a_2(1 - 2\e)}{\sqrt{\frac{\gamma}{\lambda}} + 2a_2\e}.\]
Hence
 \begin{align}
 \label{2.11}
  \begin{split} J_1 + {A_2}^\frac{2(r+1)}{p+1} J_2
  \leq \frac{1}{A_1}
\int_\Omega (u+1)^q(v+1)^p+{A_2}^\frac{2r+1-p}{p+1}\int_\Omega (u+1)^\theta
(v+1)^r + C_\e.
  \end{split}
 \end{align}
By \eqref{2.6},
$$q - \frac{\theta - 1}{2} = q + 1 - \frac{\theta + 1}{2} = (q + 1)\left[1 - \frac{p+1}{2(r+1)}\right].$$
Using Young's inequality, there holds
\begin{align*}
 & \frac{1}{A_1}\int_\Omega (u+1)^q(v+1)^p \\
= & \; \int_\Omega (u+1)^{\frac{\theta-1}{2}}(v+1)^{\frac{p-1}{2}}(u+1)^{(q + 1)\left(1 - \frac{p+1}{2(r+1)}\right)}\frac{(v+1)^{\frac{p+1}{2}}}{A_1}\\
\leq & \; \int_\Omega (u+1)^{\frac{\theta-1}{2}}(v+1)^{\frac{p-1}{2}}\left[\frac{2r+1-p}{2(r+1)}(u+1)^{q+1} + \frac{p+1}{2(r+1)} A_1^{-\frac{2(r+1)}{p+1}}(v+1)^{r+1}\right]\\
= & \; \frac{2r+1-p}{2(r+1)}J_1+\frac{p+1}{2(r+1)} A_1^{-\frac{2(r+1)}{p+1}} J_2.
\end{align*}
Similarly we have
 \begin{equation}\nonumber
  {A_2}^{\frac{2r+1-p}{p+1}}\int_\Omega (u+1)^\theta (v+1)^r \leq \frac{p+1}{2(r+1)} J_1
  + \frac{2r+1-p}{2(r+1)}{A_2}^\frac{2(r+1)}{p+1}J_2.
   \end{equation}
   Combining the above two estimates with \eqref{2.11}, we derive that
\begin{align*}
{A_2}^{\frac{2(r+1)}{p+1}}J_2\leq
\left[\frac{2r+1-p}{2(r+1)}{A_2}^{\frac{2(r+1)}{p+1}}+\frac{p+1}{2(r+1)}{A_1}^{\frac{-2(r+1)}{p+1}}\right]J_2 + C_\e,
\end{align*}
 hence
\begin{equation*}
\frac{p+1}{2(r+1)}\left[(A_1A_2)^{\frac{2(r+1)}{p+1}}-1\right]
J_2\leq C_\e.
\end{equation*}
Thus $J_2 \leq C_\e$ if $A_1A_2 > 1$. Suppose that $a_1a_2>1,$ we can take $\e > 0$ sufficiently small so that $A_1A_2>1.$

\medskip
Denote $s=r+1$. Using \eqref{2.6}, we can check directly that $a_{1}a_2>1$ is equivalent to $L(s)<0$. We conclude then for all $s > p+1$ verifying $L(s) < 0$,
there is $C_s > 0$ such that for any $(\lambda, \gamma) \in \Gamma_\sigma$,
\begin{equation}\label{2.12}
  \int_\Omega (u+1)^{\frac{\theta-1}{2}}(v+1)^{\frac{p+1}{2}}(v+1)^{s-1} = J_2 \le C_s.
\end{equation}
So we are done. \qed

\begin{rem}\label{r.2.1} Let $L$ be given by \eqref{2.1} and $H$ be given by \eqref{1.5}. A direct computation yields
$$H(x)=\left(\frac{\theta+1}{p\theta-1}\right)^4L(s), \quad \mbox{if } x=\frac{\theta+1}{p\theta-1}s.$$
Denote $s_0$ the largest root of $L$, then $x_0=\frac{\theta+1}{p\theta-1}s_0$ is the largest root of $H$, and $H(x)<0$ if and only if $L(s)<0$.
Moreover, there holds
$$L(2t_0)=\frac{16p\theta(p+1)(\theta-p)}{(\theta+1)^2}(1-2t_0).$$
So $L(2t_0)<0$ for $p < \theta$. As $\lim_{s\rightarrow\infty}L(s)= \infty,$ it follows that $2t_0 < s_0.$ If $p=\theta$, we have
\begin{align*}
L(s) = s^4 -16p^2s^2 + 32p^2s-16p^2 =(s^2+4ps-4p)(s^2-4ps +4p).
\end{align*}
For any $p > 1$, we check readily that $2t_0 = 2p+ 2\sqrt{p^2-p}$ is the largest root of $L$.
\end{rem}

\noindent
{\bf Proof of Theorem \ref{main1} completed.} Let $ (\lambda^*,\gamma^*) \in \Upsilon$ and $\sigma=\frac{\gamma^*}{\lambda^*}$. Denote
$(u,v)$ the minimal solution of \eqref{1.1} with $(\lambda, \gamma) \in \Gamma_\sigma$. Applying Lemma \ref{f}, if $p+1 <s<s_0$, there exists $C_s > 0$
such that
\[ \frac{1}{\gamma} \int_\Omega | \nabla v|^2 = \int_\Omega (u+1)^\theta v \le\int_\Omega (u+1)^\theta (v+1)^{s-1} \le C_s,\]
passing to the limit, we see that $ v^* \in H_0^1(\Omega)$. Moreover,
\[ \frac{-\Delta v^*}{\gamma^*} = (u^*+1)^\theta = \frac{ (u^*+1)^\theta}{v^*+1} v^* + \frac{ (u^*+1)^\theta}{v^*+1} \quad \mbox{in $ \Omega$.}\]
By standard elliptic theory, to show the boundedness of
$ v^*$, it is sufficient to prove that $\frac{(u^*+1)^\theta}{v^*+1} \in L^T(\Omega)$ for some $T> \frac{N}{2}$. Using (\ref{new1})
and passing to the limit, we see that there is some $ C>0$ such that
\[ \frac{ (u^*+1)^\theta}{v^*+1} \le C (u^*+1)^\frac{p \theta-1}{p+1} \quad \mbox{a.e.~in $ \Omega$.} \]  According to the estimate \eqref{2.2}
which holds also with $u^{\star},$ it follows that $ \frac{ (u^*+1)^\theta}{v^*+1} \in L^T(\Omega)$ for some $T>\frac{N}{2}$ provided
\[ \frac{(p \theta -1)}{p+1} \frac{N}{2} < \theta + \frac{( \theta +1) (s_0-1)}{p+1}, \;\; \mbox{or equivalently } N < 2 + 2 x_0
\;\mbox{where } x_0=\frac{\theta+1}{p\theta-1}s_0,\]
This is just the desired result. Moreover, using Remark \ref{r.2.1} and adopting the proof of Remark 2 in \cite{co1}, we can easily show that
\begin{align*}
x_0 \geq 2t_0\frac{\theta+1}{p\theta-1}>4,\quad \forall\,\theta\geq p>1.
  \end{align*}
This means that if $N\leq 10$, $(u^{\star},v^{\star})$ is bounded. \qed

\section{Proof of Theorem \ref{main2}}
\subsection{Some preparations}
We establish first some properties of the polynomial $L$ defined by \eqref{2.1}. Recall that without loss of generality, we can assume $1 < p \leq \theta$.
\begin{lem}\label{l.4.0}
Let $1< p \leq \theta$, then $L(2) < 0$ and $L$ has a unique root $s_0$ in $(2, \infty)$. Moreover, we have $p+1<2\theta \frac{p+1}{\theta+1}<s_0.$
\end{lem}

\noindent{\bf Proof.} As $1< p \leq \theta$, we have
\begin{align*}
L(2) & = 16 -\frac{64p\theta(p+1)}{(\theta+1)} +\frac{32p\theta(p+1)(p+\theta+2)}{(\theta+1)^2} -\frac{16p\theta(p+1)^2}{(\theta+1)^2}\\
& = 16 -\frac{64p\theta(p+1)}{(\theta+1)} + \frac{32p\theta(p+1)}{(\theta+1)} + \frac{32p\theta(p+1)^2}{(\theta+1)^2}
-\frac{16p\theta(p+1)^2}{(\theta+1)^2}\\
& = 16 -\frac{32p\theta(p+1)}{(\theta+1)} + \frac{16p\theta(p+1)^2}{(\theta+1)^2}\\
& \leq 16 - \frac{32p\theta(p+1)}{(\theta+1)} + \frac{16p\theta(p+1)}{(\theta+1)}\\
& = 16 \frac{(1 - p^2)\theta + (1 - p\theta)}{(\theta+1)} < 0.
\end{align*}
Very similarly, we can check that
\begin{align*}
L'(2) \leq 32 - \frac{32p\theta(p+1)}{(\theta+1)} < 0 \quad \mbox{and} \quad L''(s)=12s^2-\frac{32p\theta(p+1)}{\theta+1}.
\end{align*}
Then
$L''$ could change at most once the sign from negative to positive in $[2, \infty)$, hence $L$ admits a unique root in $(2,\infty).$
In addition, direct calculations yield to
\begin{align*}
L(p+1)=-\dfrac{(p+1)^2(5p\theta+\theta+p+1)(3p\theta-\theta -p-1)}{(\theta+1)^2}<0,
\end{align*}
and
\begin{align*}
-\frac{(\theta+1)^4}{16\theta(p+1)^2}L\left(2\theta \frac{p+1}{\theta+1}\right) = (3p^2-1)\theta^3 +(2p^2-p)\theta^2 -2(p^2+p)\theta +p =:K(p, \theta).
\end{align*}
It's not difficult to check that for any $1<p\leq \theta$, $K(p, \theta) > K(p, p) > 0$, hence
$L(2\theta \frac{p+1}{\theta+1})<0,$ which means $2\theta \frac{p+1}{\theta+1} < s_0$. \qed

\bigskip
Using similar ideas as in the proof for Lemma \ref{f} and following the proof of Lemma 3.1 in \cite{hhm}, we can claim
\begin{lem}
\label{l.4.1} Let $(u,v)$ be a stable solution of \eqref{1.1}. Then, for any $s > \frac{p+1}{2}$ verifying
$L(s) < 0$, there exists $C <\infty$ such that if $B_R(y) \subset \Omega$ for $R >0$, then
\begin{equation*}
\int_{B_{R/2}(y)}(u+1)^{\theta} (v+1)^{s-1}dx  \leq\frac{C}{R^2}\int_{B_R(y)}(v+1)^{s}dx.
\end{equation*}
\end{lem}
We will need also the following well known elliptic estimate. Denote $B_r := B_r(0)$ for any $r > 0$.
\begin{lem} \label{l.4.2} For $1\leq t<\frac{N}{N-2},$ there exists $C>0$ such that for any $w \in W^{2,1}(B_{2R})$ with $R > 0$, we have
\begin{equation*}
\left(\int_{B_{R}}|w|^{t} dx\right)^\frac{1}{t}\leq
CR^{N\big(\frac{1}{t}-1\big)+2}\int_{B_{2R}}|\D w|dx+CR^{N\big(\frac{1}{t}-1\big)}\int_{B_{2R}}|w|dx.
\end{equation*}
\end{lem}
As a consequence of the two above Lemmas, we state
\begin{lem}
\label{l.4.3}
 Let $(u,v)$ be a stable solution of \eqref{1.1}. Then, for any $2\leq s<\frac{N}{N-2}s_0$,
there are $\ell\in \N$ and $C > 0$ such that if $B_R(y) \subset \Omega$, then
 \begin{align*}
 \left(R^{-N}\int_{B_{2^{-\ell} R}(y)}(v+1)^{s} dx\right)^\frac{2}{s}\leq
CR^{-N}\int_{B_R(y)}(v+1)^{2} dx.
\end{align*}
\end{lem}
\textbf{Proof.} The proof of Lemma \ref{l.4.3} is very similar to that for Proposition 1 in \cite{co}, (see also Lemma 3.3 in \cite{hhm}
for a more general setting). It follows from the application of Lemma \ref{l.4.2} with $w=(v+1)^s$. We use also Lemma \ref{l.4.1} to control the integral
$$\int_{B_{R/2}(y)}(u+1)^{\theta} (v+1)^{s-1}dx$$ appeared after multiplying the equation of $v$ by $(v+1)^{s-1}\phi^2,$ where $\phi$
is a suitable cut off function. We omit the details here.\qed
\begin{rem}\label{r.4.1}
Let $ (\lambda^*,\gamma^*) \in \Upsilon$ and $\sigma=\frac{\gamma^*}{\lambda^*}$. Suppose that $(u,v)$ is a stable solution of \eqref{1.1} with $(\lambda,\gamma) \in \Gamma_\sigma.$ Although the constant C appearing in Lemma \ref{l.4.1} as well as in Lemma \ref{l.4.3} depends on $\lambda,\gamma,$ it remains bounded as $\lambda\nearrow \lambda^*.$
\end{rem}

\subsection{$\varepsilon$-regularity.} Inspired by \cite{DG}, we prove the following
$\varepsilon$-regularity result which is crucial in proving Theorem \ref{main2}. Denote
\begin{equation*}
  \alpha= \frac{2(p+1)}{p
  \theta-1},\quad \beta=\frac{2(\theta+1)}{p\theta-1},
\end{equation*}
the scaling exponents of system \eqref{1.1}.
\begin{prop}\label{p.4.1} Assume that $N\geq 2+2x_0$ and $\theta \ge p > 1$. Let $(u^{\star},v^{\star})$ be an extremal solution associated to \eqref{1.1}. There exists $\varepsilon_0>0$ such that if for some $B_{R_0}(x)\subset \Omega$ with $R_0 > 0$ and
\begin{align*}
 R_0^{2\beta-N}\int_{B_{R_0}(x)}(v^{\star}+1)^{2}\leq \varepsilon_0,
\end{align*}
then $x$ is a regular point of $(u^*, v^*)$, i.e.~$u^{\star}$, $v^{\star}$ are smooth in a neighborhood of $x$.
\end{prop}
For the proof of Proposition \ref{p.4.1}, we need to establish the following lemma.
 \begin{lem}\label{l.4.4}
  There exist $\varepsilon_1$ and $\tau \in (0, 1)$ depending on $N, p, \theta$ such that
  if $(u; v)$ is a stable solution of \eqref{1.1}, $B_{R_0}(z)\subset \Omega$ and
\begin{equation}\label{3.1}
  G_0 := R_0^{2\beta-N}\int_{B_{R_0}(z)}(v+1)^{2} dx \leq \varepsilon_1,
\end{equation}
then
\begin{equation}\label{new2}
(\tau R_0)^{2\beta-N}\int_{B_{\tau R_0}(z)}(v+1)^{2} dx\leq \frac{1}{2}G_0.
\end{equation}
\end{lem}
\textbf{Proof.} By shifting coordinates, we can assume that $z=0.$ Up to the scaling
\begin{equation}\label{3.2}
  \tilde{u}(x)+1=R_0^{\alpha}(u(R_0x)+1),\quad  \tilde{v}(x)+1=R_0^{\beta}(v(R_0x)+1),
\end{equation}
we can assume $R_0 = 1$ without loss of generality. By Lemma \ref{l.4.0}, we have $2\theta\frac{p+1}{\theta+1}<s_0$, hence by Lemma \ref{l.4.3} and \eqref{new1}, there exist $\ell\in \mathbb{N}$ and $C > 0$ such that
\begin{align*}
\int_{B_{2^{-\ell}}}(u+1)^{2\theta}\leq C\int_{B_{2^{-\ell}}}(v+1)^{2\theta\frac{p+1}{\theta+1}}
\leq C\left[\int_{B_1}(v+1)^{2}\right]^{\frac{\theta(p+1)}{\theta+1}}.
\end{align*}
Denote $r_0:= 2^{-\ell}$ and using \eqref{3.1}, we deduce that
\begin{equation}\label{3.3}
   \int_{B_{r_0}}(u+1)^{2\theta}\leq CG_0^{\frac{\theta(p+1)}{\theta+1}} \leq C\varepsilon_1^{\frac{\theta(p+1)}{\theta+1}}.
\end{equation}
Consider now the decomposition $v+1=v_1+v_2$ where
$$\left\{\begin{array}{ll}
     -\Delta v_1= 0 & \mbox{in}\;B_{r_0}\\
     v_1=v+1 & \mbox{on} \;\partial B_{r_0},
     \end{array}
     \right.
     \quad \left\{\begin{array}{ll}
   -\Delta v_2=\gamma (u+1)^\theta & \mbox{in}\;B_{r_0}\\
   v_2=0 & \mbox{on} \;\partial B_{r_0}.
\end{array}\right.$$
 Let $ 0<\tau<r_0$ (to be fixed later on), we have
\begin{equation}\label{3.4}
\int_{B_\tau}(v+1)^{2} dx\leq 2 \int_{B_\tau}v_1^{2} dx+ 2\int_{B_\tau}v_2^{2} dx.
\end{equation}
Noting that $v_1^{2}$ is subharmonic in $B_{r_0}$, we get
\begin{align}\label{3.5}
\tau ^{2\beta-N}\int_{B_\tau}v_1^{2} dx \leq C\tau ^{2\beta}\int_{B_{r_0}}v_1^{2}dx \leq C\tau ^{2\beta}\int_{B_1}(v+1)^{2} dx = C\tau^{2\beta}G_0.
\end{align}
On the other hand, by elliptic theory and \eqref{3.3}, there holds, as $G_0 \leq \varepsilon_1$,
\begin{align}\label{3.6}
  \int_{B_{r_0}} v_2^2\leq \|v_2\|_{H^2(B_{r_0})}^2 \leq CG_0^{\frac{\theta(p+1)}{\theta+1}} \leq C\varepsilon_1^{\frac{p\theta-1}{\theta+1}}G_0.
\end{align}
  Combining \eqref{3.4}-\eqref{3.6}, we obtain
   \begin{equation*}
    \tau ^{2\beta-N}\int_{B_\tau}(v+1)^{2} dx\leq C\tau^{2\beta} G_0 + C\tau^{2\beta-N}\varepsilon_1^{\frac{p\theta -1}{\theta+1}}G_0.
  \end{equation*}
  Fix $\tau >0$ so that $C \tau^{2\beta}\leq \frac14$. Then, take $\varepsilon_1>0$ sufficiently small so that $C \tau^{2\beta-N}\varepsilon_1^{\frac{p\theta-1}{\theta+1}}\leq \frac14$, we are done. \qed

\medskip\noindent
\textbf{Proof of Proposition \ref{p.4.1}.}
By approximating the extremal solution $(u^{\star}, v^{\star})$ of \eqref{1.1} by minimal solutions with parameters $(\lambda,\gamma) \in \Gamma_\sigma$, Lemma \ref{l.4.4} holds true for $v^{\star}$. As above, we can assume that $x=0$ and $R_0=1$.

\medskip
Since $N\geq 2+2x_0=2+\beta s_0$ and $s_0>2,$ we get $N-2\beta>0. $ Let $\varepsilon_1$ be the constant in Lemma \ref{l.4.4} and choose $\varepsilon_0$ such that $2^{N-2\beta}\varepsilon_0=\varepsilon_1.$ For any $y\in B_{\frac{1}{2}}$, we have
\begin{align*}
G_1 := \left(\frac{1}{2}\right)^{2\beta-N}\int_{B_{ \frac{1}{2}}(y)}(v^{\star}+1)^{2} dx\leq 2^{N-2\beta}\int_{B_1}(v^{\star}+1)^{2} dx\leq  2^{N-2\beta}\varepsilon_0 = \varepsilon_1.
\end{align*}
 Applying inductively Lemma \ref{l.4.4}, then for any $k\geq 1,$
\begin{align}\left(\frac{\tau^k}{2}\right)^{2\beta-N}\int_{B_{ \frac{\tau^k}{2}}(y)}(v^{\star}+1)^{2} dx\leq 2^{-k}G_1.
\end{align}
Let $0<\rho\leq \tau/2,$ we can take $k\in \N^*$ such that $\frac{\tau^{k+1}}{2}<\rho \leq\frac{\tau^k}{2}.$ Therefore,
\begin{align*}
\rho^{2\beta-N}\int_{B_{\rho}(y)}(v^{\star}+1)^{2} dx & \leq \left(\frac{\tau^{k+1}}{2}\right)^{2\beta-N}\int_{B_{ \frac{\tau^k}{2}}(y)}(v^{\star}+1)^{2}dx\\
& \leq 2^{N-2\beta}\tau^{2\beta-N}2^{-k}G_1 \\
& \leq C(N,\beta, \tau, \e_1)2^{-k-1}\\
& \leq C\tau^{(k+1)\delta},
\end{align*}
where $\delta = \frac{-\ln 2}{\ln \tau} > 0$. This implies that
\begin{equation}\label{3.7}
  \int_{B_{\rho}(y)}(v^{\star}+1)^{2} dx\leq C\rho^{N-2\beta+\delta}, \quad \forall\; y\in B_{\frac{1}{2}}, \; 0 < \rho \leq \frac{\tau}{2}.
\end{equation}

Furthermore, applying Lemma \ref{l.4.3} with $s=p+1$ and using \eqref{3.7}, we get
 \begin{align*}
 \int_{B_{ \rho}(y)}(v^{\star}+1)^{p+1} dx\leq  C\rho^{N-(p+1)\beta+\delta}, \quad \forall\; y\in B_{\frac{1}{2}}, \; 0 < \rho \leq \frac{\tau}{2^{\ell+1}}
\end{align*}
for some integer $\ell\geq 1.$ By approximation argument, the estimate \eqref{new1} holds a.e.~in $\Omega$, if we replace $(u, v)$ by $(u^{\star}, v^\star)$. Therefore,
 \begin{align*}\int_{B_{ \rho}(y)}(u^{\star}+1)^{\theta+1} dx\leq  C\rho^{N-(p+1)\beta+\delta}, \quad \forall\; y\in B_{\frac{1}{2}}, \; 0 < \rho \leq \frac{\tau}{2^{\ell+1}}.
\end{align*}
This means that $u^\star+1$ belongs to the Morrey space
$$L^{\theta+1,N-(p+1)\beta+\delta}(B_{\frac 12})\subset L^{\theta,N-\theta \alpha+\frac{\theta}{\theta+1}\delta}(B_{\frac 12}).$$ Finally,
applying Theorem 3.4 in \cite{BM}, we get the claim. \qed

\subsection {Proof of Theorem \ref{main2} completed}
Let $p+1<s<\frac{N}{N-2}s_0$ such that $L(s) < 0$. Let $z\in \Omega$ verifying
\begin{equation*}
  \lim_{R\rightarrow 0}R^{\beta s-N}\int_{B_{R}(z)}(v^\star+1)^sdx=0.
\end{equation*}
By H\"older's inequality, there holds
\begin{equation*}
  \lim_{R\rightarrow 0}R^{2\beta -N}\int_{B_{R}(z)}(v^\star+1)^2dx=0.
\end{equation*}
Applying Proposition \ref{p.4.1}, $z$ is a regular point for $(u^*, v^*)$. This implies that the singular set
\begin{equation*}
  \mathcal{S}\subset \left\{x\in \Omega:\limsup_{R\rightarrow 0}R^{\beta s-N}\int_{B_R(x)}(v^\star+1)^sdx>0\right\}.
\end{equation*}
Take first $s_1 > p+1$ such that $L(s_1) < 0$. Using \eqref{2.2}, for minimal solution $(u, v)$ of \eqref{1.1} with $(\lambda, \gamma) \in \Gamma_\sigma$, as $\frac{p+2s_1 - 1}{2} > 2$,
\begin{align*}\int_\Omega (v+1)^2 \leq \int_\Omega (u+1)^\frac{\theta-1}{2}(v+1)^\frac{p+2s_1-1}{2} \leq C_{s_1}.
\end{align*}
Passing to the limit, the above estimate yields that $v^* + 1 \in L^2(\Omega)$. By Lemma \ref{l.4.3}, $(v^\star+1)^s\in L_{loc}^1(\Omega)$, it follows that $\mathcal{H}^{N-\beta s}({\mathcal S}) = 0$ whenever $N - \beta s > 0$, see Theorem 5.3.4 in \cite{D}. Tending $s$ to $\frac{N}{N-2}s_0, $ we conclude that the Hausdorff dimension of $\mathcal{S}$ is less or equal than $\max(N-\frac{2N}{N-2}x_0, 0)$, recalling that $x_0= \frac{\theta+1}{p\theta-1}s_0$. As $N\geq 2+2x_0,$ the claim follows.\qed

\bigskip\noindent
{\bf Acknowledgments}. I would like to thank Professor Dong Ye for suggesting me this problem and for many helpful comments.

\end{document}